\documentclass{amsart}
\usepackage{amsmath}
\usepackage{amsthm}
\usepackage{amsfonts}
\usepackage{amssymb}
\newtheorem{Theorem}{Theorem}[section]

\newtheorem{Lemma}[Theorem]{Lemma}

\theoremstyle{definition}

\numberwithin{equation}{section}

\begin{document}
\title[Discrete and embedded eigenvalues]{Discrete and embedded
eigenvalues for one-dimensional Schr\"odinger operators}

\author{Christian Remling}

\address{Mathematics Department\\
University of Oklahoma\\
Norman, OK 73019-0315}
\email{cremling@math.ou.edu}
\urladdr{www.math.ou.edu/$\sim$cremling}
\date{\today}

\begin{abstract}
I present an example of a discrete Schr\"odinger operator
that shows that it is possible to have embedded singular spectrum and, at the same time,
discrete eigenvalues that approach the edges of the essential spectrum (much) faster
than exponentially. This settles a conjecture of Simon (in the negative).
The potential is of von Neumann-Wigner type, with careful navigation around
a previously identified borderline situation.
\end{abstract}
\maketitle
\section{Introduction}
I am interested in one-dimensional discrete Schr\"odinger equations,
\begin{equation}
\label{se}
u(n+1)+u(n-1)+V(n)u(n)=Eu(n),
\end{equation}
and the associated self-adjoint operators
\[
(Hu)(n) = \begin{cases} u(n+1)+u(n-1)+V(n)u(n) & \quad (n\ge 2) \\
u(2)+V(1)u(1) & \quad (n=1) \end{cases}
\]
on $\ell_2(\mathbb N)$. We could also consider whole line operators
(on $\ell_2(\mathbb Z)$), and for the purposes of this paper,
that would actually make very little difference.

Recent work has shown that there are fascinating and unexpected relations between the
discrete and essential spectrum of $H$. If $V\equiv 0$, then $\sigma_{ac}(H) =[-2,2]$,
$\sigma_{sing}(H)=\emptyset$. It turns out that perturbations
that substantially change the character of the spectrum on $[-2,2]$ must also
introduce new spectrum outside this interval. Indeed,
Damanik and Killip \cite{DK} proved the spectacular result that if $\sigma\setminus
[-2,2]$ is finite, then $[-2,2]\subset\sigma$ and the spectrum continues to be
purely absolutely continuous on $[-2,2]$. The situation when
$\sigma\setminus [-2,2]$ is a possibly infinite set of discrete eigenvalues,
with $\pm 2$ being the only possible accumulation points, was subsequently
investigated by Damanik and myself \cite{DR} (honesty demands that I point out that
we actually treated the
analogous problems in the continuous setting).

A natural question, which was not addressed in \cite{DK,DR}, concerns
the minimal assumptions that will still imply that the spectrum
is purely absolutely continuous on $[-2,2]$ in these situations.
Put differently:
\textit{How fast can the discrete eigenvalues
approach the edges of the essential spectrum $\sigma_{ess}=[-2,2]$ if the essential spectrum
is} not \textit{purely absolutely continuous? Is there in fact any bound on this
rate of convergence?}

So we assume that $\sigma\setminus [-2,2]=\{ E_n \}$, and we introduce
\[
d_n \equiv \textrm{dist }(E_n, [-2,2]) .
\]
We also assume that $\{E_n\}$ is infinite and $d_n\to 0$.
It would in fact be natural (but not really necessary) to arrange
the eigenvalues so that $d_1\ge d_2\ge \ldots$.

The hunt for examples with some singular spectrum on $[-2,2]$, but
rapidly decreasing $d_n$'s was opened by Damanik, Killip, and Simon in \cite{DKS}.
This paper has an example
where $d_n \lesssim e^{-cn}$ and $0\in\sigma_{pp}$. Based on this,
Simon conjectured \cite[Sect.\ 13.5]{SimOPUC} that the slightly stronger condition
$\sum 1/\ln d_n^{-1} <\infty$ might suffice to conclude that the spectrum is purely absolutely
continuous on $[-2,2]$. The purpose of this paper is to improve on the
counterexample from \cite{DKS}; this will also show that the above
conjecture needs to be adjusted. More precisely, we will prove:
\begin{Theorem}
\label{T1.1}
There exists a potential $V$ so that:
(1) For $E=0$, the Schr\"odinger equation \eqref{se} has an $\ell_2$ solution $u$.\\
(2) $d_n \le e^{-cn^2}$ for some $c>0$.
\end{Theorem}

Such a potential $V$ is in fact explicitly given by
\begin{equation}
\label{defV}
V(n)=\frac{(-1)^n}{n} \left( 1 + \frac{2}{\ln n} \right) \quad\quad (n\ge 3) ;
\end{equation}
here $2$ could be replaced by any other constant $c>1$. If a complete
definition of $V$ is desired, we can put $V(1)=V(2)=0$. However, the behavior of $V$ on finite sets
is quite irrelevant for what we do here. Note also in this context that by adjusting $V(1)$, say,
we can achieve that $0\in\sigma_{pp}$.

To motivate \eqref{defV}, let us for a moment consider the simpler potential
$V(n)=g(-1)^n/n$. This is basically a discrete variant of the classical von Neumann-Wigner
potential \cite{vNW}. The values $g=\pm 1$ for the coupling constant are critical in two senses:
First of all, there exists a square summable solution to the Schr\"odinger equation \eqref{se}
at energy $E=0$ if and only if $|g|>1$. Second, if $|g|\le 1$, then the operator $H$ has no spectrum
outside $[-2,2]$ \cite[Proposition 5.9]{DHKS}. On the other hand, if $|g|>1$, there must be infinitely many eigenvalues
$E_n$ with $|E_n|>2$ by the result from \cite{DK} discussed above. A rather detailed
analysis is possible, and $V(n)=g(-1)^n/n$ is in fact the example of Damanik-Killip-Simon
mentioned above: The eigenvalues approach $\pm 2$ exponentially fast \cite[Theorem 1]{DKS}.

So it seems to make sense to make $g$ $n$-dependent and approach the threshold value $g=1$ more cautiously.
The aim of this paper is to show that this idea works. Of course, since there are no positive
results beyond the Damanik-Killip theorem at this point, the question of what the fastest possible
decay of the $d_n$'s is must remain open. The fact that the type of counterexample used here feels right
together with some experimentation with the $V^2/4$ trick from \cite{DHKS} have in fact led me to believe
that the rate $d_n\lesssim e^{-cn^2}$ might already be the correct answer, but this is probably too
bold a claim.

The plan of this paper is as follows: We will prove the estimate on the eigenvalues (part (2) of
Theorem \ref{T1.1}) in Sect.\ 2--4. We will use oscillation theory: Roughly speaking, it is possible
to locate eigenvalues by counting zeros. Our basic strategy is in part inspired by the
treatment of \cite{KS}. From a more technical point of view, the statement we formulate as
Lemma \ref{L3.1} is very much at the heart of the matter.
Part (1) of Theorem \ref{T1.1} will be proved in Sect.\ 5. We will use
a discrete version of Levinson's Theorem (compare \cite[Theorem 1.3.1]{East}) as our main tool.
\section{Variation of constants}
We will write $V(n)=V_0(n)+V_1(n)$ with $V_0(n)=(-1)^n/n$ and $E=2+\epsilon$
and view $V_1$ as well as $\epsilon$ as perturbations. This strategy seems
especially appropriate here because one can in fact solve \eqref{se} with $V=V_0$
and $E=2$ (almost) explicitly. This observation, which is crucial for what follows,
is from \cite{DHKS}. There is actually no need to discuss the equation for a general
$E\ge 2$ here. Rather, it suffices to have good control on the solutions for $E=2$
because we can then refer to oscillation theory at a later stage of the proof.

Let us begin with the unperturbed problem: So, consider \eqref{se} with
$V(n)=V_0(n)=(-1)^n/n$ and $E=2$. As observed in \cite{DHKS}, if we define
\begin{equation}
\label{2.2}
\varphi_{2n}=\varphi_{2n+1} = \prod_{j=1}^n \left( 1+ \frac{1}{2j-1} \right) ,
\end{equation}
then $\varphi_n$ solves this equation. We find a second, linearly independent solution
$\psi$ to the same equation by using constancy of the Wronskian,
\begin{equation}
\label{2.1}
\varphi_n \psi_{n+1} - \varphi_{n+1} \psi_n = 1 ,
\end{equation}
and making the ansatz $\psi_n=C_n\varphi_n$. Plugging this into \eqref{2.1}, we see that
$\psi$ will solve \eqref{se} if
\begin{equation}
\label{2.9}
C_{n+1}-C_n = \frac{1}{\varphi_n\varphi_{n+1}} .
\end{equation}
For later use, we record asymptotic formulae for these solutions.
A warning may be in order here: What we call $\varphi$ in Lemma \ref{L2.1} below
differs from the $\varphi$ defined in \eqref{2.2} by a constant factor. By the
same token, the asymptotic formula for $C_n$ of course implies a particular choice
of the constant in the general solution of \eqref{2.9}.
\begin{Lemma}
\label{L2.1}
There exist solutions $\varphi$, $\psi_n=C_n\varphi_n$ to \eqref{se} with
$V=V_0$ and $E=2$ satisfying the following asymptotic formulae:
\[
\varphi_{2n}=\varphi_{2n+1} = (2n)^{1/2} + O(n^{-1/2}),\quad C_n = \ln n + O(1/n)
\]
\end{Lemma}
\begin{proof}[Sketch of proof]
Take logarithms in \eqref{2.2} and asymptotically evaluate the resulting
sum by using Taylor expansions and approximating sums by
integrals. Then use this information to
analyze \eqref{2.9}.
\end{proof}
To analyze the full equation, with $V$ given by \eqref{defV}, we use variation
of constants. So write $T_0(n)$ for the transfer matrix of the unperturbed problem,
that is,
\[
T_0(n) = \begin{pmatrix} \varphi_n & \psi_n \\ \varphi_{n+1} & \psi_{n+1} \end{pmatrix} .
\]
Then $\det T_0(n)=1$ because of \eqref{2.1}, and
\[
T_0^{-1}(n) = \begin{pmatrix} \psi_{n+1} & -\psi_n \\ -\varphi_{n+1} & \varphi_n \end{pmatrix} .
\]
It will also be convenient to write $V(n)=V_0(n)+(-1)^n W_n$, so
\begin{equation}
\label{2.5}
W_n = \frac{2}{n\ln n} .
\end{equation}

Now let $y$ be a solution of the Schr\"odinger equation \eqref{se} with $E=2$ and potential \eqref{defV},
and introduce $D_n\in\mathbb C^2$ by writing
\[
\begin{pmatrix} y_n \\ y_{n+1} \end{pmatrix} = T_0(n) D_n .
\]
A calculation then shows that $D_n$ solves
\begin{equation}
\label{2.3}
D_n-D_{n-1} = A_nD_{n-1},\quad\quad A_n \equiv (-1)^nW_n\varphi_n^2 \begin{pmatrix}
C_n & C_n^2 \\ -1 & - C_n \end{pmatrix} .
\end{equation}
We have used the fact that $\psi_n=C_n\varphi_n$.

We notice that $A_{n+1}\approx -A_n$, so we expect at least partial
cancellations. To exploit this, we perform \textit{two} steps (rather than just one)
in the iteration \eqref{2.3}. Clearly, $D_{2n+1}=(1+A_{2n+1})(1+A_{2n})D_{2n-1}$,
so we define
\[
M_n =(1+A_{2n+1})(1+A_{2n}) .
\]
Using the formulae $\varphi_{2n+1}=\varphi_{2n}$ and $C_{2n+1}=C_{2n}+\varphi_{2n}^{-2}$,
we find that
\begin{multline*}
M_n = 1 + (W_{2n}-W_{2n+1}+W_{2n}W_{2n+1})\varphi_{2n}^2
\begin{pmatrix} C_{2n} & C_{2n}^2 \\ -1 & - C_{2n} \end{pmatrix}\\
-W_{2n+1} \begin{pmatrix} 1 & 2C_{2n} + \varphi_{2n}^{-2} \\ 0 & -1 \end{pmatrix}
+W_{2n}W_{2n+1} \begin{pmatrix} 1 & C_{2n} \\ 0 & 0 \end{pmatrix} .
\end{multline*}
We see at this point already that the idea of doing two steps at once was a major success
because this new matrix $M_n$ differs from the unity matrix by a correction of order
$O(\ln n/n)$ whereas $A_n$ itself only satisfies $\|A_n\| = O(\ln n)$.
For the following calculations, it will be convenient to introduce some abbreviations and write $M_n$
in the form
\begin{equation}
\label{eqM}
M_n = \begin{pmatrix} 1+\epsilon_n - w_n + \rho_n & c_n(\epsilon_n - 2w_n +\rho_n-\rho'_n) \\
-\epsilon_n/c_n & 1-\epsilon_n + w_n \end{pmatrix} ,
\end{equation}
where
\begin{gather}
\label{2.6}
\epsilon_n = (W_{2n}-W_{2n+1} + W_{2n}W_{2n+1} ) \varphi^2_{2n} C_{2n},\quad\quad
w_n = W_{2n+1},\\
\label{2.7}
c_n=C_{2n},\quad\quad
\rho_n = W_{2n}W_{2n+1}, \quad\quad
\rho'_n= \frac{W_{2n+1}}{\varphi_{2n}^2 C_{2n}} .
\end{gather}
\section{Counting zeros}
We now want to use the difference equations from the preceding section to derive
upper bounds on the number of zeros (more precisely: sign changes) of the solution
$y$ on large intervals. Our goal in this section is to prove the following:
\begin{Lemma}
\label{L3.1}
There exist $n_0\in\mathbb N$ and $A>0$ so that the following holds: If $N_1, N_2\in
\mathbb N$ with $N_2\ge N_1\ge n_0$ and
\[
\ln N_2 \le \ln N_1 + A \ln^{1/2} N_1,
\]
then there exists a solution $y$ of \eqref{se} with $E=2$ and
potential \eqref{defV} satisfying $y_n>0$ for $N_1\le n \le N_2$.
\end{Lemma}
We start out by finding the eigenvalues and eigenvectors of $M_n$ from \eqref{eqM}.
The eigenvalues are given by
\[
\lambda_{\pm}(n) = 1 + \frac{\rho_n}{2} \pm w_n \sqrt{ 1 + \frac{\epsilon_n\rho'_n}{w_n^2}
- \frac{\rho_n}{w_n} + \frac{\rho_n^2}{4w_n^2}} .
\]
We will need information on the asymptotic behavior. From Lemma \ref{L2.1} and the
definitions (see \eqref{2.5}, \eqref{2.6}, and \eqref{2.7}), we obtain that
\begin{equation}
\label{asymp}
\epsilon_n = \frac{1}{n} + O\left( \frac{1}{n\ln n} \right), \quad\quad
\rho_n, \rho'_n = O\left( \frac{1}{n^2\ln^2 n} \right) .
\end{equation}
This shows, first of all, that
\begin{equation}
\label{asympev}
\lambda_{\pm}(n) = 1 \pm \frac{1}{n\ln 2n} + O\left( \frac{1}{n^2\ln n} \right) .
\end{equation}
Next, write the eigenvector corresponding to $\lambda_+$ in the form
$v_+ = \bigl( \begin{smallmatrix} 1 \\ -a_n \end{smallmatrix} \bigr)$. It then follows
from \eqref{eqM} that $a_n$ satisfies
\[
-\frac{\epsilon_n}{c_n} + (1-\epsilon_n+w_n-\lambda_+(n))(-a_n)=0 ,
\]
and by using the slightly more precise formula
\[
\lambda_+(n)= 1 + w_n + \frac{\epsilon_n\rho'_n}{2w_n} + O\left( \frac{1}{n^2\ln^2 n} \right)
\]
instead of \eqref{asympev}, we obtain from this that
\begin{equation}
\label{3.1}
a_n = \frac{1}{c_n} - \frac{1}{4n\ln^2 n} + O\left( \frac{1}{n\ln^3 n} \right) .
\end{equation}

We are interested in the sign of $y_n$, so obviously only the \textit{direction} of $D_n$ matters.
Let $\theta_n$ be the angle that $D_{2n-1}$ makes with the eigenvector $v_+$; see also
Figure 1 below.\\[0.8cm]

\setlength{\unitlength}{1cm}
\begin{picture}(10,8)
\put(0,8){Figure 1}
\put(1,3.6){\line(1,0){8}}
\put(5,0.5){\line(0,1){6.5}}
{\thicklines
\put(5,3.6){\vector(4,-1){3}}
\put(5,3.6){\vector(2,3){2}}}
\linethickness{0.1mm}
\bezier{300}(6.39,5.68)(7.8,4.5)(7.43,2.99)
\put(7.2,6.5){$D_{2n-1}$}
\put(8.2,2.7){$v_+(n)$}
\put(6.2,4){$\theta_n$}
\end{picture}\\[0.6cm]

\begin{Lemma}
\label{L3.2}
There exists $n_0\in\mathbb N$ such that the following holds: If $n\ge n_0$ and
$0\le\theta_n\le \pi/2$, then $y_{2n-1}>0$ and $y_{2n}>0$.
\end{Lemma}
\begin{proof}
The condition on $\theta_n$ implies that we can write
$D_{2n-1}= k_n\bigl( \begin{smallmatrix} 1 \\ -d_n \end{smallmatrix} \bigr)$
with $k_n>0$ and $d_n\le a_n$. Now
\begin{align*}
\begin{pmatrix} y_{2n-1} \\ y_{2n} \end{pmatrix} & = T_0(2n-1)D_{2n-1}
= k_n \begin{pmatrix} \varphi_{2n-2} & \varphi_{2n-2} C_{2n-1} \\
\varphi_{2n} & \varphi_{2n}C_{2n} \end{pmatrix} \begin{pmatrix} 1 \\ -d_n \end{pmatrix} \\
& = k_n \begin{pmatrix} \varphi_{2n-2}(1-C_{2n-1}d_n) \\ \varphi_{2n} (1-C_{2n}d_n) \end{pmatrix}
\end{align*}
Since $k_n,\varphi_{2n-2},\varphi_{2n}>0$ and $C_{2n}>C_{2n-1}>0$ (we may have to take $n$
sufficiently large here), we see that $y_{2n-1}$, $y_{2n}$ will certainly
be positive if $1-C_{2n}d_n=1-c_nd_n>0$ or, equivalently, $1/c_n>d_n$. But \eqref{3.1} shows that $1/c_n>a_n$
for large $n$, and, as noted above, $d_n\le a_n$, so this condition holds.
\end{proof}
To motivate the subsequent arguments, we now make some preliminary, informal remarks about
the dynamics of the recursion $D_{2n+1}=M_n D_{2n-1}$:
First of all, a calculation shows that the eigenvector $v_-=v_-(n)$ associated with the
small eigenvalue $\lambda_-<1$ is of the form $v_-=\bigl( \begin{smallmatrix} 1 \\ -b_n \end{smallmatrix}
\bigr)$ with $b_n>a_n$, so $v_-$ lies below $v_+$. However, $b_n-a_n = O(\ln^{-2}n)$, so $v_+$
and $v_-$ are almost parallel.

Now an application of a $2\times 2$ matrix moves the
vector towards the eigenvector corresponding to the large eigenvalue. So in the case at hand,
$D$ will approach $v_+$, or, in other words, $\theta_n$ will decrease. At the same time, $v_+(n)$
approaches the positive $x$-axis, but this is a comparatively small effect. Nevertheless, a crossing
between $D$ and $v_+$ will eventually occur. Our task is to bound from below the number of iterations
it takes (starting from $\theta=\pi/4$, say) to reach this crossing.

We will use the eigenvector $v_+=v_+(n)=\bigl( \begin{smallmatrix} 1 \\ -a_n \end{smallmatrix}
\bigr)$ and the orthogonal vector
$\bigl( \begin{smallmatrix} a_n \\ 1 \end{smallmatrix} \bigr)$ as our basis of $\mathbb R^2$.
As $\theta_n$ was defined as the
angle between $D_{2n-1}$ and $v_+(n)$, it follows that $D_{2n-1}$ is a constant multiple
of the vector
\begin{equation}
\label{3.8}
\cos\theta_n v_+(n) + \sin\theta_n \begin{pmatrix} a_n \\ 1 \end{pmatrix} .
\end{equation}
Conversely, we can find $\theta$ using the fact that $D$ has such a representation.
More precisely, to compute $\theta_{n+1}$ from $\theta_n$, we apply the matrix $M_n$ to
the vector from \eqref{3.8} and
then take scalar products with $v_+(n+1)$ and $\bigl( \begin{smallmatrix}
a_{n+1} \\ 1 \end{smallmatrix} \bigr)$. These operations produce multiples
of $\cos\theta_{n+1}$ and $\sin\theta_{n+1}$, respectively. We omit the details of this
routine calculation. The result is as follows: If we introduce $t_n=\tan\theta_n$, then
\begin{equation}
\label{eqt}
t_{n+1} = \frac{s_nt_n + \lambda_+(n) (a_{n+1}-a_n)}{\widetilde{s}_nt_n + \lambda_+(n)
(1+a_n a_{n+1})} ,
\end{equation}
where $a_n$ was defined above (see also \eqref{3.1}) and
\begin{align*}
s_n & = (a_{n+1}, 1)
M_n \begin{pmatrix} a_n \\ 1 \end{pmatrix} ,\\
\widetilde{s}_n & = (1, -a_{n+1})
M_n \begin{pmatrix} a_n \\ 1 \end{pmatrix} .
\end{align*}
From \eqref{eqM}, \eqref{asymp}, \eqref{3.1}, and Lemma \ref{L2.1},
we obtain the asymptotic formulae
\begin{equation}
\label{asymps}
s_n = 1+ a_na_{n+1} + O\left( \frac{1}{n\ln n} \right),\quad\quad
\widetilde{s}_n = \frac{\ln n}{n} + O\left( \frac{1}{n} \right) .
\end{equation}
We will prove Lemma \ref{L3.1} by analyzing the recursion \eqref{eqt}. As a preliminary,
we observe the following:
\begin{Lemma}
\label{L3.3}
There exists $n_0\in\mathbb N$ so that the following holds: If $n\ge n_0$ and $0\le t_n\le 1/\ln n$,
then also $t_{n+1}\le 1/\ln (n+1)$.
\end{Lemma}
\begin{proof}
Let
\[
f(x) = \frac{s_n x + \lambda_+ (a_{n+1}-a_n)}{\widetilde{s}_n x + \lambda_+
(1+a_n a_{n+1})}
\]
be the function from \eqref{eqt}. Then
\[
f'(x) = \frac{\lambda_+}{\left( \widetilde{s}_n x + \lambda_+
(1+a_n a_{n+1}) \right)^2} \left( s_n (1+a_na_{n+1}) - \widetilde{s}_n(a_{n+1}-a_n) \right) ,
\]
and since $a_{n+1}-a_n =O(1/(n\ln^2 n))$, the derivative is positive for large enough $n$. Therefore,
$t_{n+1} = f(t_n) \le f(1/\ln n)$, and by dividing through by $1+a_na_{n+1}$ and using
\eqref{asympev}, \eqref{3.1}, and \eqref{asymps}, we see that
\[
t_{n+1} \le \frac{\left( 1 +O\left( \frac{1}{n\ln n} \right)\right) \frac{1}{\ln n}
+O\left( \frac{1}{n\ln^2 n} \right)}{\frac{1}{n}+1+O\left( \frac{1}{n\ln n} \right)}
= \left( 1- \frac{1}{n} \right) \frac{1}{\ln n} + O\left( \frac{1}{n\ln^2 n} \right) .
\]
On the other hand,
\[
\ln (n+1) = \ln n + \ln \left( 1+\frac{1}{n} \right) = \ln n + O\left( \frac{1}{n} \right) ,
\]
so the claim follows.
\end{proof}
We are now ready for the \textit{proof of Lemma \ref{L3.1}.} We must show that when solving
the basic recursion $D_{2n+1}=M_nD_{2n-1}$ (or its variant \eqref{eqt}), at least as much time as
specified in the statement is spent in the region where $y_n>0$. We will in fact show that $t_n$
spends such an amount of time already in the region where
\begin{equation}
\label{3.9}
\ln^{-3/2} n \le t_n \le \ln^{-1} n .
\end{equation}
Lemma \ref{L3.2} shows that this condition indeed implies that $y_{2n-1},y_{2n}>0$.
In fact, \eqref{3.9} might look unnecessarily restrictive so that the whole analysis
would appear to be rather crude.
However, an argument similar to the one we are
about to give shows that the time spent in the neglected regions is at most of the same order
of magnitude. More precisely, if $0\le t_n \le \ln^{-3/2} n$ for $N_1 \le n\le N_2$,
then these $N_1$, $N_2$ also satisfy the estimate from Lemma \ref{L3.1}. Moreover, if $\ln^{-1} n\le
t_n\le M$ for $N_1\le n\le N_2$, then $N_2 \le CN_1$, with $C$ independent of $M>0$.
These remarks, together with a more careful analysis of the crossing between
$D$ and $v_+$, show that the condition from Lemma \ref{L3.1} is sharp.

Let us now proceed with the strategy outlined above. Assume that \eqref{3.9} holds.
From \eqref{eqt}, we obtain that
\[
t_{n+1}-t_n = \frac{\left[ s_n - \lambda_+(n)(1+a_na_{n+1}) \right] t_n
+\lambda_+(n) (a_{n+1}-a_n) - \widetilde{s}_n t_n^2}{\widetilde{s}_nt_n + \lambda_+(n) (1+a_na_{n+1})} .
\]
Note that $s_n - \lambda_+(n)(1+a_na_{n+1})=O(1/(n\ln n))$. Also, $t_n\le\ln^{-1} n$ by assumption,
so the first term in the numerator is of the order $O(1/(n\ln^2 n))$. As for the next term, recall that
$a_{n+1}-a_n =O(1/(n\ln^2 n))$. Finally, since we are assuming that
$t_n\ge \ln^{-3/2}n$, we have that $\widetilde{s}_n t_n^2 \gtrsim 1/(n\ln^2 n)$,
so up to a constant factor, this term is not smaller than the other two summands from the numerator.
The denominator clearly is of the form $1+o(1)$. So, putting things together, we see that
\[
t_{n+1}-t_n \ge -C \frac{\ln n}{n} t_n^2,
\]
with $C>0$. It will in fact be convenient to write this in the form
\[
t_{n+1}-t_n \ge -C \frac{\ln n}{n} t_n t_{n+1} ,
\]
with an adjusted constant $C>0$. We can then introduce $r_n=1/t_n$, and this new
variable obeys
\begin{equation}
\label{diffr}
r_{n+1}-r_n \le C \frac{\ln n}{n} .
\end{equation}
This was derived under the assumptions that $n$ is sufficiently large and that we have the two bounds
\[
\ln n\le r_n\le \ln^{3/2} n .
\]
Our final task is to use \eqref{diffr} to
find an estimate on the first $n$ for which the second inequality fails to hold.
Recall also that Lemma \ref{L3.3} says that there can't be any such problems
with the first inequality.
Suppose that $N_1\in\mathbb N$ is sufficiently large and $r_{N_1}=\ln N_1$.
We can proceed by induction: What we have just shown
says that if $r_j\le \ln^{3/2} j$ for $j=N_1,N_1+1,\ldots , n-1$, then
$r_n\ge \ln n$ and
\begin{equation}
\label{3.5}
r_n \le \ln N_1 + C \sum_{j=N_1}^{n-1} \frac{\ln j}{j} .
\end{equation}
So we can keep going as long as the right-hand side of
\eqref{3.5} is $\le \ln^{3/2} n$. Now clearly
\[
\sum_{j=N_1}^{n-1} \frac{\ln j}{j} \lesssim \ln^2 n - \ln^2 N_1 ,
\]
so, recalling from the statement of Lemma \ref{L3.1} what we are actually trying to
prove, we see that it suffices to show that given a constant $C>0$, we can find $A>0$, $n_0\in\mathbb N$
so that if $N_1\ge n_0$, then the condition
\begin{equation}
\label{3.6}
\ln N_2 \le \ln N_1 + A\ln^{1/2} N_1
\end{equation}
implies that
\begin{equation}
\label{3.7}
\ln N_1 + C(\ln^2 N_2 -\ln^2 N_1) \le \ln^{3/2} N_2 .
\end{equation}
Indeed, it will then follow that $y_n>0$ for $n=2N_1-1,\ldots , 2N_2$, and a
simple adjustment gives Lemma \ref{L3.1} as originally stated, without the
factors of $2$.

But the above claim is actually quite obvious: By taking squares in \eqref{3.6}, we obtain
the estimate
\[
\ln^2 N_2 - \ln^2 N_1 \le 2A \ln^{3/2}N_1 + A^2 \ln N_1 \le 3A \ln^{3/2} N_1
\]
(say), and if also $3AC<1$, then \eqref{3.7} follows at once.
\hfill$\square$
\section{Oscillation theory}
In this section, we will use Lemma \ref{L3.1} to derive the desired estimate
on the discrete eigenvalues. For the time being, we are concerned with
eigenvalues $E_n>2$; in particular, we then have that $d_n=E_n-2$. We will need
some standard facts from oscillation theory; for proofs, we refer the reader
to \cite[Chapter 4]{Te}. This reference gives a careful discussion of all
the results we will need (and several others), but some caution is required when looking up results
in \cite{Te} because the operators discussed there are the negatives of the operators
generated by \eqref{se}.

We first state a couple of comparison theorems: If $u_1$, $u_2\not\equiv 0$ both solve the same
Schr\"odinger equation \eqref{se}, then the number of zeros (more precisely:
sign changes) on any fixed interval differs by at most one \cite[Lemma 4.4]{Te}.
Also, if $E\le E'$
and $u$, $u'$ solve the Schr\"odinger equation with energies $E$ and $E'$, respectively, and these
solutions have the same initial phase at $N_1$ (i.e.\ $u(N_1+1)/u(N_1)=u'(N_1+1)/u'(N_1)$),
then for any $N_2>N_1$, $u$ has at least as many zeros on $\{N_1,\ldots , N_2 \}$ as $u'$.
This can be deduced from \cite[Theorem 4.7]{Te}.

Finally, the following is a more direct consequence of \cite[Theorem 4.7]{Te}:
If $u$ solves \eqref{se}, $u(0)=0$, $u(1)=1$, and $u$ has
$N$ zeros on $\mathbb N$, then $E_N\le E< E_{N-1}$. Here, we of course implicitly
assume that the eigenvalues $E_n$ are arranged
in their natural order: $E_0>E_1>E_2>\ldots$ (and $E_{-1}:=\infty$).
In other words, it is possible to locate discrete eigenvalues by counting zeros.

Let us now use these facts to derive Theorem \ref{T1.1}(2) from Lemma \ref{L3.1}.
Define $N_n=\exp(\gamma n^2)+x_n$, with $0<\gamma<A^2/4$, where $A$ is the constant from
Lemma \ref{L3.1}; $x_n\in [0,1)$ is chosen
so that $N_n\in\mathbb N$. Then $\ln N_{n+1}\le \ln N_n
+ A \ln^{1/2}N_n$ for all sufficiently large $n$, so Lemma \ref{L3.1} and the facts
reviewed in the preceding paragraphs imply that any nontrivial solution $u$ of \eqref{se}
with $E=2$ has at most $C+n$ zeros on $\{1,\ldots, N_n \}$, where $C$ is a fixed
constant, independent of $n$. Moreover, the same holds for any nontrivial solution
of \eqref{se} with $E\ge 2$ because increasing the energy leads to fewer zeros by
the comparison theorem quoted above.

Now fix $E>2$ and let $d=E-2$; for convenience, also assume that $d<1/2$, say. Then if $m\ge 3/d$,
then $|V(m')|<d$ for all $m'\ge m$, and thus the operator on $\ell_2(\{ m,m+1,\ldots \} )$ doesn't have
any spectrum in $[2+d,\infty)$. As a consequence, any solution $u$ of the Schr\"odinger equation
\eqref{se} with $E=2+d$ has at most one zero on $\{ m,m+1,\ldots \}$. If this is combined
with what has been observed above, it follows that an arbitrary non-trivial solution $u$
to \eqref{se} with $E=2+d$ has at most $C + n(d)$ zeros on
$\mathbb N$; here $n(d)$ must be chosen so that $N_{n(d)}\ge 3/d$.
In other words, it is possible to pick $n(d)\lesssim \ln^{1/2}d^{-1}$.

To sum this up: Any non-trivial solution $u$ to \eqref{se} with $E=2+d$ has at most
$C_1 + C_2\ln^{1/2} d^{-1} \le C_3 \ln^{1/2} d^{-1}$ zeros on $\mathbb N$. We can now use
that part of oscillation theory that relates the location of the eigenvalues to the number of zeros.
It follows that if $N=[C\ln^{1/2} d^{-1}]$, then the $N$th eigenvalue, $E_N$, satisfies $E_N\le 2+d$.
By rearranging, we see that $E_N-2\le \exp(-cN^2)$, as claimed in part (2) of Theorem \ref{T1.1}.

Of course, we haven't talked about eigenvalues $<-2$ yet, but this part of the claim is established
by a completely analogous analysis. We can use the fact that $-E$ is an eigenvalue of
the original problem if and only if
$E$ is an eigenvalue for the potential $-V$. This reduces matters to the discussion
of the eigenvalues bigger than $2$, but for the two potentials $V$ and $-V$. We have just discussed $V$,
and, not surprisingly, it turns out that the sign change is quite irrelevant and we can simply
run the whole argument again, with only a few very minor modifications.
So we will not give any details, and these general remarks conclude the proof of
part (2) of Theorem \ref{T1.1}.
\section{Asymptotic integration}
What we will do below is modelled on the treatment of similar problems in
the continuous setting. See in particular \cite[Sect.\ 4.3]{East}.

We want to analyze the solutions of the discrete Schr\"odinger equation \eqref{se}
with potential \eqref{defV} and $E=0$. For the following computations, it will be
convenient to write $V(n)=(-1)^n 2v_n$, so
\[
v_n = \frac{1}{2n} \left( 1+ \frac{2}{\ln n} \right)
\]
for $n\ge 3$. Given a solution $y$ of
\[
y_{n+1}+y_{n-1}+(-1)^n 2v_ny_n=0 ,
\]
introduce $Y_n=\bigl( \begin{smallmatrix} y_{n-1} \\ y_n \end{smallmatrix} \bigr)$.
Then $Y$ solves
\begin{equation}
\label{5.2}
Y_{n+1} = \begin{pmatrix} 0 & 1 \\ -1 & (-1)^{n+1} 2v_n \end{pmatrix} Y_n .
\end{equation}
We again use a variation of constants type transformation, treating $V$ as the
perturbation. So define a new variable $Z$ by $Y_n=T_n Z_n$, where
\[
T_n = \begin{pmatrix} \cos \pi (n-1)/2 & \sin \pi (n-1)/2 \\ \cos \pi n/2 & \sin \pi n/2
\end{pmatrix} ;
\]
note that indeed $T_{n+1} = \bigl( \begin{smallmatrix} 0 & 1 \\ -1 & 0 \end{smallmatrix}
\bigr) T_n$, that is, $T$ solves the unperturbed equation. A calculation shows that $Z$ obeys
\[
Z_{n+1} = Z_n - v_n \begin{pmatrix} 0 & 1+ (-1)^{n+1} \\ 1+(-1)^n & 0
\end{pmatrix} Z_n .
\]
Here, we have used the fact that the trigonometric functions
only take the values $0$, $\pm 1$ at integer multiples
of $\pi/2$, and, for example, $\cos^2 n\pi/2 =(1+(-1)^n)/2$.

Next, we diagonalize the non-oscillating part of the perturbation: To this end, write
\[
Z_n = \begin{pmatrix} 1 & 1 \\ 1 & -1 \end{pmatrix} W_n .
\]
Then $W$ solves
\[
W_{n+1}=W_n + v_n \begin{pmatrix} -1 & (-1)^{n+1} \\ (-1)^n & 1 \end{pmatrix} W_n .
\]
Finally, we can now approximately get rid of the oscillating part with the help of a
transformation of the type
\[
W_n = (1+v_nA_n) U_n .
\]
The matrix $A_n$ will be chosen shortly; it will satisfy $A_n =O(1)$. Since $v_n^2,
v_{n+1}-v_n \in\ell_1$, we then have that
\[
(1+v_{n+1}A_{n+1})^{-1} = 1 -v_n A_{n+1} + B_n
\]
with $B_n\in\ell_1$. It thus follows that
\[
U_{n+1} = U_n + v_n\left[ A_n-A_{n+1} +
\begin{pmatrix} -1 & (-1)^{n+1} \\ (-1)^n & 1 \end{pmatrix} \right] U_n + R_n U_n ,
\]
with $R_n\in\ell_1$. This suggests that we take
\[
A_n = \frac{(-1)^n}{2} \begin{pmatrix} 0 & 1 \\ -1 & 0 \end{pmatrix} ,
\]
and this choice leads to the equation
\begin{equation}
\label{eqU}
U_{n+1} = \begin{pmatrix} 1 - v_n & 0 \\ 0 & 1 + v_n \end{pmatrix} U_n + R_nU_n
\end{equation}
for $U$.

We expect that the summable perturbation $R$ does not change the asymptotics of the solutions.
We are especially interested in the decaying solution, and we want to show that, more specifically,
there exists a solution $U$ satisfying $U_n\approx \left( \prod_{j=j_0}^{n-1} (1 -v_j)
\right) \bigl( \begin{smallmatrix} 1 \\ 0 \end{smallmatrix} \bigr)$.
To do this, we mimic the proof of Levinson's Theorem. We will basically follow the
presentation given in \cite{BR} (but see also \cite[Sect.\ 1.4]{East}). When appropriate,
we will right away specialize to the case at hand although the underlying arguments
are actually of a much more general character throughout.

Consider the ``integral equation''
\begin{equation}
\label{eqC}
C_n = \begin{pmatrix} 1 \\ 0 \end{pmatrix} - \sum_{j=n}^{\infty}
\frac{1}{1-v_j}\begin{pmatrix} 1 & 0 \\ 0 & \prod_{k=n}^j \frac{1-v_k}{1+v_k} \end{pmatrix} R_j C_j
\equiv e_1 - (TC)_n
\end{equation}
We will see later that this is basically a way of rewriting \eqref{eqU}.
Since $R\in\ell_1$, $0<v_n\le c<1$, and $v_n\to 0$, the sum from \eqref{eqC} defines a bounded operator on
$\ell_{\infty}(\{ j_0, j_0+1, \ldots \} ;\mathbb C^2)$ (bounded sequences on
$\{n: n\ge j_0 \}$ taking values in $\mathbb C^2$). More precisely,
\[
\|TC\|_{\infty} \le \frac{1}{1-v_{j_0}}
\sum_{j=j_0}^{\infty} |R_j| \cdot \|C\|_{\infty} .
\]
So if we take $j_0$ sufficiently large, then in fact $\|T\|<1$ and thus
$1+T$ is boundedly invertible on $\ell_{\infty}(\{ j_0, j_0+1, \ldots \} ;\mathbb C^2 )$.
In particular, $C\equiv (1+T)^{-1}e_1$
is a bounded solution to \eqref{eqC}. This boundedness of course also makes sure that the series
from \eqref{eqC} converges.

Moreover, it follows from \eqref{eqC} that $C_n\to e_1$ as $n\to\infty$. Finally, as already
announced, we obtain a solution to the original equation \eqref{eqU} from this $C$: define
\begin{equation}
\label{5.1}
U_n = \left[ \prod_{j=j_0}^{n-1} (1-v_j) \right] C_n \equiv p_n C_n.
\end{equation}
Then $U_n$ solves \eqref{eqU} for $n\ge j_0$. To verify this claim, call the
diagonal matrix from \eqref{eqU} $\Lambda_n$, so that \eqref{eqU} becomes
$U_{n+1}=(\Lambda_n+R_n)U_n$. Next observe that $p_{n+1}e_1 = p_n \Lambda_n e_1$ and
\[
p_{n+1} \begin{pmatrix} 1 & 0 \\ 0 & \prod_{k=n+1}^j \frac{1-v_k}{1+v_k} \end{pmatrix}
= p_n \Lambda_n \begin{pmatrix} 1 & 0 \\ 0 & \prod_{k=n}^j \frac{1-v_k}{1+v_k} \end{pmatrix} .
\]
So if $U_n$ is defined by \eqref{5.1} and $C_n \in\ell_{\infty}$ solves \eqref{eqC}, then
\[
U_n = p_n e_1 - p_n \sum_{j=n}^{\infty} \frac{1}{1-v_j}
\begin{pmatrix} 1 & 0 \\ 0 & \prod_{k=n}^j \frac{1-v_k}{1+v_k} \end{pmatrix} R_j C_j ,
\]
thus
\begin{align*}
U_{n+1} & = p_{n+1} e_1 - p_{n+1} \sum_{j=n+1}^{\infty} \frac{1}{1-v_j}
\begin{pmatrix} 1 & 0 \\ 0 & \prod_{k=n+1}^j \frac{1-v_k}{1+v_k} \end{pmatrix} R_j C_j \\
& = \Lambda_n p_n e_1 - p_n \Lambda_n \sum_{j=n+1}^{\infty} \frac{1}{1-v_j}
\begin{pmatrix} 1 & 0 \\ 0 & \prod_{k=n}^j \frac{1-v_k}{1+v_k} \end{pmatrix} R_j C_j \\
& = \Lambda_n p_n e_1 - p_n \Lambda_n \sum_{j=n}^{\infty} \frac{1}{1-v_j}
\begin{pmatrix} 1 & 0 \\ 0 & \prod_{k=n}^j \frac{1-v_k}{1+v_k} \end{pmatrix} R_j C_j\\
& \quad
+ p_n \Lambda_n \frac{1}{1-v_n}
\begin{pmatrix} 1 & 0 \\ 0 & \frac{1-v_n}{1+v_n} \end{pmatrix} R_n C_n\\
& = \Lambda_n U_n + p_n R_n C_n = (\Lambda_n +R_n) U_n ,
\end{align*}
as required.

We can now go back to the original variable $Y$:
\[
Y_n = T_n \begin{pmatrix} 1 & 1 \\ 1 & -1 \end{pmatrix} (1+v_n A_n) U_n
\]
See also \eqref{5.2}. Since $U_n = p_n (e_1 + o(1))$ as $n\to\infty$, it follows that
there is a solution $y$ of the Schr\"odinger equation satisfying
\[
|y_n| = (1+o(1)) \prod_{j=j_0}^n (1-v_j) \quad\quad\quad (n\to\infty) .
\]
Now
\begin{align*}
\ln \prod_{j=j_0}^n (1-v_j) & = \sum_{j=j_0}^n \ln (1-v_j) = -\sum_{j=j_0}^n v_j +O(1)\\
& = -\sum_{j=j_0}^n \left( \frac{1}{2j} + \frac{1}{j\ln j} \right) +O(1)
= -\frac{1}{2}\ln n - \ln\ln n +O(1),
\end{align*}
so
\[
y_n^2 \lesssim \frac{1}{n\ln^2 n} ,
\]
and $y$ is indeed square summable. The proof of Theorem \ref{T1.1} is complete.
\hfill$\square$

\end{document}